\documentclass[final,3p,times,nopreprintline]{elsarticle}

\usepackage{amsmath,amssymb,amsthm,mathtools}
\usepackage{microtype}
\allowdisplaybreaks[2]
\newtheorem{theorem}{Theorem}

\newcommand{\pint}{\int_\Omega\frac{|u|^p}{|x|^{s_1}}\,dx}
\newcommand{\qint}{\int_\Omega\frac{|u|^q}{|x|^{s_2}}\,dx}

\begin{document}

\begin{frontmatter}

\title{Large-parameter nonexistence for all finite-energy solutions of a boundary Hardy--Sobolev problem}

\author[aff1]{Zhi-Yun Tang}
\ead{tzy@swu.edu.cn}

\author[aff1]{Ying Lv\corref{cor1}}
\ead{ly0904@swu.edu.cn}

\address[aff1]{School of Mathematics and Statistics, Southwest University, Chongqing, 400715, People's Republic of China}

\cortext[cor1]{Corresponding author.}

\begin{abstract}
Let $\Omega\subset\mathbb R^N$ ($N\geq3$) be a bounded domain with $0\in\partial\Omega$. We study a semilinear Dirichlet equation with a critical Hardy--Sobolev source $|x|^{-s_2}|u|^{q-2}u$ and an absorption term $-\lambda|x|^{-s_1}|u|^{p-2}u$, where $0\leq s_1<s_2<2$, $q=2^*(s_2)<p<2^*(s_1)$ and $2^*(s)=2(N-s)/(N-2)$. If $\partial\Omega$ is of class $C^{1,\alpha}$ near $0$, $0<\alpha\leq1$, and $p\geq2^*(s_1)-2\alpha/(N-2)$, then for all sufficiently large $\lambda>0$ the only solution in $H_0^1(\Omega)$ is zero. No curvature condition or boundary regularity away from $0$ is required. The proof combines a localized Pohozaev inequality with a weighted energy estimate based on the angular Poincar\'e inequality.
\end{abstract}


\end{frontmatter}

\section{Introduction}

Let $N\geq3$, let $\Omega\subset\mathbb R^N$ be a bounded domain with $0\in\partial\Omega$, and set $2^*(s)=2(N-s)/(N-2)$ for $0\leq s<2$. Fix $0\leq s_1<s_2<2$, put $q=2^*(s_2)$, and consider
\begin{equation}\label{eq:main}
\left\{
\begin{aligned}
-\Delta u&=\frac{|u|^{q-2}u}{|x|^{s_2}}-\lambda\frac{|u|^{p-2}u}{|x|^{s_1}}&&\text{in }\Omega,\\
u&=0&&\text{on }\partial\Omega,
\end{aligned}
\right.
\qquad u\in H_0^1(\Omega),
\end{equation}
where $\lambda>0$ and
\begin{equation}\label{eq:exponents}
q=2^*(s_2)<p<2^*(s_1).
\end{equation}

Li and Lin \cite{Li2012} raised the positive-solution problem for two competing Hardy--Sobolev terms. On bounded $C^1$ domains with $C^2$ boundary near $0$ and negative mean curvature there, Cerami, Zhong, and Zou \cite{Cerami2015} obtained positive solutions for small $\lambda>0$ when $q<p\leq2^*(s_1)$, and for almost every $\lambda>0$ when $q<p<2^*(s_1)-2/(N-2)$.

In the doubly critical case $p=2^*(s_1)$, Wang and Su \cite[Corollary~1.5]{Wang2025} obtained an explicit upper bound on $\lambda$ for nonzero solutions on smooth domains with negative mean curvature at $0$ and in the half-space. Their quasilinear results also include small-parameter existence in the radial whole-space setting, the half-space and open cones. Tang and Tang \cite{Tang2025} proved large-parameter nonexistence of nonzero solutions at the doubly critical endpoint and studied the continuation of positive very weak solutions for certain singly critical exponents. Yang \cite{Yang2026} relaxed the exponent restriction in the latter result using boundary Poisson barriers.

For subcritical source exponents, Tang and Tang \cite{Tang2025a} obtained at least two positive solutions for small $\lambda>0$ and nonexistence of nonzero solutions for large $\lambda$ when $2<p\leq2^*(s_1)$ and $2<q<2+(2-s_2)(p-2)/(2-s_1)$. In a whole-space variant with an additional linear term $+u$ on the left-hand side, they proved nonexistence of nonzero solutions for every $\lambda>0$ when $q=2^*(s_2)<p\leq2^*(s_1)$ and obtained positive solutions under suitable subcritical exponent conditions \cite{Tang2025c}.

For equations with a single critical Hardy--Sobolev term, the boundary geometry affects the best constant and its attainability \cite{Ghoussoub2004,Ghoussoub2006}; local concavity can lead to multiplicity, whereas star-shapedness with respect to the singular point can rule out nontrivial solutions \cite{Ghoussoub2006a}. For the two-term equation on bounded domains star-shaped with respect to $0$, Tang and Tang \cite{Tang2025b} also obtained positive solutions in a suitable subcritical range.

In the singly critical range \eqref{eq:exponents}, Xu, Zhang, and Zhong \cite{Xu2027} gave an explicit sufficient small-parameter criterion for positive solutions on $C^1$ domains with $C^2$ boundary near $0$, provided that the Hardy--Sobolev best constant of the domain is smaller than that of the half-space. For $C^{1,\alpha}$ boundary near $0$ with $0<\alpha\leq1$, He and Lv \cite[Theorem~1.1(i)]{He2026} proved existence of positive solutions for every $\lambda>0$ when $q<p<2^*(s_1)-2\alpha/(N-2)$, assuming a homogeneous leading boundary term of degree $1+\alpha$ and negative $\alpha$-mean curvature.

Within this singly critical range, Tang and Tang \cite[Theorem~1.5]{Tang2025b} proved nonexistence of nonzero solutions for every $\lambda>0$ on bounded domains star-shaped with respect to $0$. He and Lv \cite[Theorem~1.1(ii)]{He2026} showed that, under the same local $C^{1,\alpha}$ boundary condition and exponent threshold, the set of parameters admitting positive solutions is bounded.

In this paper, we prove that for sufficiently large $\lambda$ the only solution in $H_0^1(\Omega)$ is zero, so the conclusion covers all finite-energy solutions, not only positive ones. The distinction from Tang and Tang is the trade-off between global star-shapedness and the large-parameter restriction; the distinction from He and Lv is the solution class. Our proof combines a localized Pohozaev inequality with a weighted energy estimate based on the angular Poincar\'e inequality, without relying on the uniform pointwise a priori estimate or concentration analysis used by He and Lv.

\begin{theorem}\label{thm:main}
Assume \eqref{eq:exponents} and suppose that $\partial\Omega$ is of class $C^{1,\alpha}$ in a neighborhood of $0$, where $0<\alpha\leq1$. If
\[
p\geq2^*(s_1)-\frac{2\alpha}{N-2},
\]
then there exists $\Lambda>0$, depending only on $N,s_1,s_2,p,\alpha$ and $\Omega$, such that for every $\lambda>\Lambda$, the only solution of \eqref{eq:main} in $H_0^1(\Omega)$ is $u\equiv0$.
\end{theorem}

No curvature sign or boundary regularity away from $0$ is required.

\section{Proof of Theorem~\ref{thm:main}}

\begin{proof}
Let $u\in H_0^1(\Omega)$ be a nonzero real-valued solution of \eqref{eq:main}. Throughout the proof, $C>0$ denotes a constant independent of $u$ and $\lambda$, which may change from line to line. Put
\[
\theta:=N-s_1-\frac{N-2}{2}p
=\frac{N-2}{2}\bigl(2^*(s_1)-p\bigr).
\]
The assumptions give $0<\theta\leq\alpha\leq1$. By the Hardy--Sobolev inequality and H\"older's inequality, the weighted integrals below are finite. Testing \eqref{eq:main} with $u$ gives the identity
\begin{equation}\label{eq:nehari}
\int_\Omega|\nabla u|^2\,dx+\lambda \pint=\qint.
\end{equation}

\medskip
\noindent\textit{Step 1. A localized Pohozaev inequality.}
After a rotation, take a neighborhood $U$ of $0$ such that
\[
\Omega\cap U=\{x=(x',x_N)\in U:x_N>h(x')\},
\qquad h\in C^{1,\alpha},\quad h(0)=0,\quad \nabla h(0)=0.
\]
Let $\chi\in C_c^\infty(U)$ be nonnegative and equal to $1$ near $0$. For a sufficiently large fixed $M>0$, define
\[
Y(x)=\chi(x)\bigl(x-M|x|^{1+\alpha}e_N\bigr).
\]
On the graph $x_N=h(x')$, writing $\nu=(\nabla h,-1)/\sqrt{1+|\nabla h|^2}$ for the outer unit normal, we have
\[
Y\cdot\nu
=\frac{\chi(x)\bigl(x'\cdot\nabla h(x')-h(x')+M|x|^{1+\alpha}\bigr)}
{\sqrt{1+|\nabla h(x')|^2}}\geq0.
\]
Indeed, the $C^{1,\alpha}$ assumption gives $|h(x')|\leq C|x'|^{1+\alpha}$ and $|\nabla h(x')|\leq C|x'|^\alpha$, so a fixed $M$ suffices. Writing $DY$ for the Jacobian of $Y$ and $I_N$ for the identity matrix, we have
\begin{equation}\label{eq:Yestimate}
|DY-I_N|+|\operatorname{div}Y-N|
+\left|\frac{Y\cdot x}{|x|^2}-1\right|
\leq C|x|^\alpha.
\end{equation}
Near $0$ this follows by differentiating $x-M|x|^{1+\alpha}e_N$. Away from $0$, the estimate follows after increasing $C$, since $Y$ is supported in a fixed compact subset of the boundary chart.

For the localized Pohozaev identity (cf.\ \cite{Pucci1986}), we first establish the required regularity. By Kato's inequality \cite[Corollary~1.2]{Brezis2004},
\[
-\Delta|u|\leq |x|^{-s_2}|u|^{q-2}|u|
\]
in the weak sense. Away from $0$, locally in $\Omega$ and up to the regular boundary chart, the coefficient $|x|^{-s_2}|u|^{q-2}$ belongs to $L^{N/(2-s_2)}$, because
\[
(q-2)\frac{N}{2-s_2}=\frac{2N}{N-2}=2^*(0),
\qquad \frac{N}{2-s_2}>\frac N2.
\]
Since $|u|\in H_0^1(\Omega)$, the local boundedness estimates \cite[Theorems~8.17 and~8.25]{Gilbarg2001} give local boundedness of $u$ in $\Omega$ and up to the graph away from $0$. The right-hand side of \eqref{eq:main} is thus locally bounded there. By \cite[Corollary~8.36]{Gilbarg2001}, $u$ is locally $C^{1,\gamma}$ in $\Omega$ and up to the graph away from $0$, for some $\gamma>0$. The right-hand side is then locally H\"older continuous in $\Omega$, so interior Poisson regularity gives $u\in C^2_{\mathrm{loc}}(\Omega)$. These regularity bounds may depend on the particular solution and are used only to justify the identity.

Let $\psi\in C^\infty([0,\infty))$ satisfy $0\leq\psi\leq1$, $\psi=0$ on $[0,1]$, and $\psi=1$ on $[2,\infty)$, and put
\[
Y_\varepsilon(x)=\psi(|x|/\varepsilon)Y(x).
\]
Consider the vector field
\[
(Y_\varepsilon\cdot\nabla u)\nabla u
-\frac12|\nabla u|^2Y_\varepsilon
+\left(\frac{|u|^q}{q|x|^{s_2}}-\frac{\lambda|u|^p}{p|x|^{s_1}}\right)Y_\varepsilon.
\]
On taking its divergence, the term
\[
(Y_\varepsilon\cdot\nabla u)\left(\Delta u+\frac{|u|^{q-2}u}{|x|^{s_2}}
-\lambda\frac{|u|^{p-2}u}{|x|^{s_1}}\right)
\]
vanishes by \eqref{eq:main}. The boundary term is supported on the regular graph away from $0$. Since $u=0$ there, $\nabla u=(\partial_\nu u)\nu$. Integrating the negative divergence over $\Omega$, we obtain
\begin{align}
&\int_\Omega\left(\frac{\operatorname{div}Y_\varepsilon}{2}|\nabla u|^2
-(DY_\varepsilon\nabla u)\cdot\nabla u\right)dx\notag\\
&\quad+\frac{\lambda}{p}\int_\Omega
\left(\operatorname{div}Y_\varepsilon-s_1\frac{Y_\varepsilon\cdot x}{|x|^2}\right)
\frac{|u|^p}{|x|^{s_1}}\,dx\notag\\
&\quad-\frac1q\int_\Omega
\left(\operatorname{div}Y_\varepsilon-s_2\frac{Y_\varepsilon\cdot x}{|x|^2}\right)
\frac{|u|^q}{|x|^{s_2}}\,dx
=-\frac12\int_{\partial\Omega\cap U}(Y_\varepsilon\cdot\nu)|\partial_\nu u|^2\,d\sigma\leq0.
\label{eq:pohozaev-eps}
\end{align}
To justify this integration, restrict to the boundary chart and first integrate above the shifted graphs $x_N=h(x')+t$ and then let $t\downarrow0$. The uniform bounds $|Y_\varepsilon|\leq C|x|$ and $|DY_\varepsilon|\leq C$ allow passage to the limit $\varepsilon\downarrow0$ in the volume integrals by dominated convergence and the finite-energy/Hardy--Sobolev bounds. Since $Y_\varepsilon\cdot\nu\geq0$, no integrability of the boundary term at $0$ is needed.

By the definitions of $q$ and $\theta$, we have
\[
\frac{N-s_2}{q}=\frac{N-2}{2},
\qquad
\frac{N-s_1}{p}=\frac{N-2}{2}+\frac{\theta}{p}.
\]
By \eqref{eq:Yestimate}, the limit of the left-hand side of \eqref{eq:pohozaev-eps} differs from
\[
\frac{N-2}{2}\int_\Omega|\nabla u|^2\,dx
+\lambda\left(\frac{N-2}{2}+\frac{\theta}{p}\right)\pint
-\frac{N-2}{2}\qint
\]
by at most
\[
C\int_\Omega|x|^\alpha
\left(|\nabla u|^2+\lambda\frac{|u|^p}{|x|^{s_1}}
+\frac{|u|^q}{|x|^{s_2}}\right)dx.
\]
The expression above equals $\frac{\lambda\theta}{p}\pint$ by \eqref{eq:nehari}. We conclude that
\begin{equation}\label{eq:pohozaev-ineq}
\frac{\lambda\theta}{p}\pint
\leq C\int_\Omega|x|^\alpha
\left(|\nabla u|^2+\lambda\frac{|u|^p}{|x|^{s_1}}
+\frac{|u|^q}{|x|^{s_2}}\right)dx.
\end{equation}

\medskip
\noindent\textit{Step 2. A weighted energy estimate.}
Write $B_R=B_R(0)$. Choose $R>0$ sufficiently small, with $R<1$, so that for every $0<r<R$ the complement of $\Omega\cap\partial B_r$ in $\partial B_r$ has surface measure at least $|\partial B_r|/4$. To see this, choose a spherical cap $\{\omega\in\mathbb S^{N-1}:\omega_N<-c\}$ of measure greater than $|\mathbb S^{N-1}|/4$. Since $h(r\omega')=O(r^{1+\alpha})$, this cap lies below the graph for all sufficiently small $r$.

Set
\[
w(x)=\min\{|x|,R\}^{\theta},
\qquad v(x)=\sqrt{w(x)}\,u(x).
\]
Since $0<\theta\leq\alpha\leq1$,
\begin{equation}\label{eq:w-properties}
|x|^\alpha\leq Cw(x),\qquad w(x)\leq|x|^\theta.
\end{equation}
Moreover, $wu,v\in H_0^1(\Omega)$: $w,\sqrt w$ are bounded and $|\nabla w|+|\nabla\sqrt w|\leq C|x|^{-1}$. Replacing $|x|$ by $\max\{|x|,\delta\}$ in the weights gives bounded Lipschitz multipliers; the resulting products converge in $H^1(\Omega)$ as $\delta\downarrow0$ by the Hardy inequality for the zero extension and dominated convergence.

Testing \eqref{eq:main} with $wu$ and using
\[
\nabla v=\sqrt w\,\nabla u+u\nabla\sqrt w,
\qquad
|\nabla\sqrt w|^2=\frac{\theta^2}{4}\frac{w}{|x|^2}
\quad\text{a.e. in }\Omega\cap B_R,
\]
we obtain
\begin{equation}\label{eq:weighted-test-identity}
\int_\Omega|\nabla v|^2\,dx
+\lambda\int_\Omega w\frac{|u|^p}{|x|^{s_1}}\,dx
=\int_\Omega w\frac{|u|^q}{|x|^{s_2}}\,dx
+\frac{\theta^2}{4}\int_{\Omega\cap B_R}\frac{v^2}{|x|^2}\,dx.
\end{equation}

Extend $v$ by zero to $\mathbb R^N$. For almost every $r\in(0,R)$, the function $f(\omega)=v(r\omega)$ belongs to $H^1(\mathbb S^{N-1})$ and vanishes a.e. on $\{\omega\in\mathbb S^{N-1}:r\omega\notin\Omega\}$, a set of measure at least $|\mathbb S^{N-1}|/4$. If $\bar f$ denotes its spherical average, then
\[
|\mathbb S^{N-1}|\,|\bar f|^2
=\frac1{|\mathbb S^{N-1}|}\left|\int_{\{\omega\in\mathbb S^{N-1}:r\omega\in\Omega\}}f\,d\sigma\right|^2
\leq\frac34\int_{\mathbb S^{N-1}}|f|^2\,d\sigma.
\]
The first nonzero eigenvalue of $-\Delta_{\mathbb S^{N-1}}$ is $N-1$ (the eigenvalues are $k(k+N-2)$, $k=0,1,\ldots$). The Poincar\'e inequality on the sphere therefore gives
\[
\int_{\mathbb S^{N-1}}|\nabla_\omega f|^2\,d\sigma
\geq(N-1)\int_{\mathbb S^{N-1}}|f-\bar f|^2\,d\sigma\geq\frac{N-1}{4}\int_{\mathbb S^{N-1}}|f|^2\,d\sigma.
\]
Integrating in polar coordinates, we find
\begin{equation}\label{eq:angular}
\frac{\theta^2}{4}\int_{\Omega\cap B_R}\frac{v^2}{|x|^2}\,dx
\leq\frac{\theta^2}{N-1}\int_\Omega|\nabla v|^2\,dx
\leq\frac12\int_\Omega|\nabla v|^2\,dx,
\end{equation}
where the last inequality follows from $\theta\leq1$ and $N\geq3$. Combining \eqref{eq:weighted-test-identity} and \eqref{eq:angular}, we obtain
\begin{equation}\label{eq:weighted-test}
\frac12\int_\Omega|\nabla v|^2\,dx
+\lambda\int_\Omega w\frac{|u|^p}{|x|^{s_1}}\,dx
\leq\int_\Omega w\frac{|u|^q}{|x|^{s_2}}\,dx.
\end{equation}

Let $t=(q-2)/(p-2)\in(0,1)$. From the definition of $\theta$ and $q=2^*(s_2)$,
\[
q=2(1-t)+pt,
\qquad s_2=2(1-t)+(s_1+\theta)t.
\]
Since $v^2=wu^2$ and $w\leq|x|^\theta$, it follows pointwise that
\[
\frac{w|u|^q}{|x|^{s_2}}
=\left(\frac{v^2}{|x|^2}\right)^{1-t}
\left(\frac{|u|^p}{|x|^{s_1}}\right)^t
\left(\frac{w}{|x|^\theta}\right)^t
\leq\left(\frac{v^2}{|x|^2}\right)^{1-t}
\left(\frac{|u|^p}{|x|^{s_1}}\right)^t.
\]
H\"older's inequality and the Hardy inequality for the zero extension of $v$ therefore yield
\begin{equation}\label{eq:holder}
\int_\Omega w\frac{|u|^q}{|x|^{s_2}}\,dx
\leq C\left(\int_\Omega|\nabla v|^2\,dx\right)^{1-t}\left(\pint\right)^t.
\end{equation}
By \eqref{eq:weighted-test} and \eqref{eq:holder},
\[
\frac12\int_\Omega|\nabla v|^2\,dx+\lambda\int_\Omega w\frac{|u|^p}{|x|^{s_1}}\,dx
\leq C\left(\int_\Omega|\nabla v|^2\,dx\right)^{1-t}\left(\pint\right)^t.
\]
Young's inequality gives
\[
\frac14\int_\Omega|\nabla v|^2\,dx+\lambda\int_\Omega w\frac{|u|^p}{|x|^{s_1}}\,dx
\leq C\pint.
\]
In particular, $\|\nabla v\|_{L^2(\Omega)}^2\leq C\int_\Omega |x|^{-s_1}|u|^p\,dx$, and \eqref{eq:holder} then also gives
\[
\int_\Omega w\frac{|u|^q}{|x|^{s_2}}\,dx\leq C\pint.
\]
Thus
\begin{equation}\label{eq:weighted-energy}
\int_\Omega|\nabla v|^2\,dx
+\lambda\int_\Omega w\frac{|u|^p}{|x|^{s_1}}\,dx
+\int_\Omega w\frac{|u|^q}{|x|^{s_2}}\,dx
\leq C\pint.
\end{equation}
Furthermore, since $\sqrt w\,\nabla u=\nabla v-u\nabla\sqrt w$, \eqref{eq:angular} gives
\[
\int_\Omega w|\nabla u|^2\,dx
\leq2\int_\Omega|\nabla v|^2\,dx+\frac{\theta^2}{2}\int_{\Omega\cap B_R}\frac{v^2}{|x|^2}\,dx
\leq3\int_\Omega|\nabla v|^2\,dx.
\]
Combining this with \eqref{eq:w-properties} and \eqref{eq:weighted-energy}, we arrive at
\begin{equation}\label{eq:weighted-final}
\int_\Omega|x|^\alpha
\left(|\nabla u|^2+\lambda\frac{|u|^p}{|x|^{s_1}}
+\frac{|u|^q}{|x|^{s_2}}\right)dx
\leq C\pint.
\end{equation}

\medskip
\noindent\textit{Step 3. Conclusion.} Combining \eqref{eq:pohozaev-ineq} and \eqref{eq:weighted-final}, we obtain
\[
\frac{\lambda\theta}{p}\pint\leq C\pint.
\]
Since $u\not\equiv0$ and $\theta>0$, division by the positive integral gives
\[
\lambda\leq\frac{Cp}{\theta}.
\]
The constant is independent of the solution and of $\lambda$. Taking $\Lambda=Cp/\theta$ proves Theorem~\ref{thm:main}.
\end{proof}

\section*{Acknowledgments}
This work was supported by the National Natural Science Foundation of China (No. 12601208), the Postdoctoral Fellowship Program of CPSF (No. GZC20261664), and the Innovation Research 2035 Pilot Plan of Southwest University (No. SWU-XDPY22015).

\section*{Declaration of generative AI and AI-assisted technologies in the manuscript preparation process}

During manuscript preparation, the authors used ChatGPT (OpenAI) as an auxiliary tool for literature searches, preliminary checks of mathematical arguments, and editing of the exposition; they reviewed and edited the AI-assisted material as needed and take full responsibility for the final manuscript.

\bibliographystyle{elsarticle-num}
\bibliography{references}

@Article{Li2012,
  author     = {Li, Y. Y. and Lin, C.-S.},
  journal    = {Arch. Ration. Mech. Anal.},
  title      = {{A Nonlinear Elliptic PDE with Two Sobolev-Hardy Critical Exponents}},
  year       = {2012},
  number     = {3},
  pages      = {943--968},
  volume     = {203},
  doi        = {10.1007/s00205-011-0467-2},
}

@Article{Cerami2015,
  author     = {Cerami, G. and Zhong, X. and Zou, W.},
  journal    = {Calc. Var. Partial Differential Equations},
  title      = {On some nonlinear elliptic {PDE}s with {S}obolev-{H}ardy critical exponents and a {L}i-{L}in open problem},
  year       = {2015},
  number     = {2},
  pages      = {1793--1829},
  volume     = {54},
  doi        = {10.1007/s00526-015-0844-z},
}

@Article{Tang2025,
  author   = {Tang, Zhi-Yun and Tang, Xianhua},
  journal  = {J. Differential Equations},
  title    = {On {L}i-{L}in's open problem},
  year     = {2025},
  pages    = {113244},
  volume   = {435},
  doi      = {10.1016/j.jde.2025.113244},
}

@Article{Ghoussoub2004,
  author     = {Ghoussoub, N. and Kang, X. S.},
  journal    = {Ann. Inst. H. Poincar\'{e} C Anal. Non Lin\'{e}aire},
  title      = {Hardy-{S}obolev critical elliptic equations with boundary singularities},
  year       = {2004},
  number     = {6},
  pages      = {767--793},
  volume     = {21},
  doi        = {10.1016/j.anihpc.2003.07.002},
}

@Article{Ghoussoub2006,
  author     = {Ghoussoub, N. and Robert, F.},
  journal    = {Geom. Funct. Anal.},
  title      = {The effect of curvature on the best constant in the {H}ardy-{S}obolev inequalities},
  year       = {2006},
  number     = {6},
  pages      = {1201--1245},
  volume     = {16},
  doi        = {10.1007/s00039-006-0579-2},
}

@Article{Ghoussoub2006a,
  author   = {Ghoussoub, N. and Robert, F.},
  journal  = {IMRP Int. Math. Res. Pap.},
  title    = {Concentration estimates for {E}mden-{F}owler equations with boundary singularities and critical growth},
  year     = {2006},
  pages    = {21867, 1--85},
  volume   = {2006},
  doi      = {10.1155/IMRP/2006/21867},
}

@article {Xu2027,
    AUTHOR = {Xu, Z. and Zhang, J. and Zhong, X.},
     TITLE = {An explicit small-parameter criterion for the {L}i--{L}in
              bad-sign problem},
   JOURNAL = {Appl. Math. Lett.},
    VOLUME = {186},
      YEAR = {2027},
     PAGES = {Paper No. 110181},
       DOI = {10.1016/j.aml.2026.110181},
}

@misc{He2026,
  author       = {He, Q. and Lv, Z.},
  title        = {A sharp threshold in the {Li--Lin} open problem},
  year         = {2026},
  howpublished = {arXiv:2611.30561v1},
}

@Book{Gilbarg2001,
  author    = {Gilbarg, D. and Trudinger, N. S.},
  publisher = {Springer-Verlag},
  title     = {Elliptic partial differential equations of second order},
  year      = {2001},
  address   = {Berlin},
  note      = {Reprint of the 1998 edition},
  series    = {Classics in Mathematics},
  pages     = {xiv+517},
}

@article {Yang2026,
    AUTHOR = {Yang, C. W.},
     TITLE = {Boundary {P}oisson barriers and a sharper continuation theorem
              for {L}i-{L}in's problem},
   JOURNAL = {Int. J. Differ. Equ.},
    VOLUME = {2026},
      YEAR = {2026},
     PAGES = {Paper No. e4318320, 12},
       DOI = {10.1155/ijde/4318320},
}

@Article{Tang2025b,
  author   = {Tang, Zhi-Yun and Tang, Xianhua},
  journal  = {Calc. Var. Partial Differential Equations},
  title    = {A positive solution of the elliptic equation on a starshaped domain with boundary singularities},
  year     = {2025},
  number   = {9},
  pages    = {Paper No. 285},
  volume   = {64},
  doi      = {10.1007/s00526-025-03162-8},
}

@Article{Pucci1986,
  author     = {Pucci, P. and Serrin, J.},
  journal    = {Indiana Univ. Math. J.},
  title      = {A general variational identity},
  year       = {1986},
  number     = {3},
  pages      = {681--703},
  volume     = {35},
  doi        = {10.1512/iumj.1986.35.35036},
}

@Article{Tang2025a,
  author   = {Zhi-Yun Tang and Xianhua Tang},
  journal  = {Commun. Anal. Mech.},
  title    = {Positive solutions for critical singular elliptic equations without {A}mbrosetti-{R}abinowitz type conditions},
  year     = {2025},
  number   = {2},
  pages    = {462--473},
  volume   = {17},
  doi      = {10.3934/cam.2025019},
}

@Misc{Tang2025c,
  author        = {Zhi-Yun Tang and Xianhua Tang},
  title         = {The {Li-Lin's open problem on $\mathbb{R}^N$}},
  year          = {2025},
  archiveprefix = {arXiv},
  eprint        = {2505.03613},
  primaryclass  = {math.AP},
}

@Article{Wang2025,
  author   = {Wang, C. and Su, J.},
  journal  = {NoDEA Nonlinear Differential Equations Appl.},
  title    = {On the double weighted critical quasilinear {H}\'enon problems},
  year     = {2025},
  number   = {2},
  pages    = {Paper No. 19, 23},
  volume   = {32},
  doi      = {10.1007/s00030-025-01028-8},
}

@article {Brezis2004,
    AUTHOR = {Br\'{e}zis, H. and Ponce, A. C.},
     TITLE = {Kato's inequality when {$\Delta u$} is a measure},
   JOURNAL = {C. R. Math. Acad. Sci. Paris},
    VOLUME = {338},
      YEAR = {2004},
    NUMBER = {8},
     PAGES = {599--604},
       DOI = {10.1016/j.crma.2003.12.032},
}

\end{document}